\documentclass[a4paper,12pt]{article}

\makeatletter
\@addtoreset{footnote}{page}
\makeatother

\usepackage{style-enumitem}

\usepackage{amsthm}
\usepackage{amsmath,amssymb,latexsym,amsfonts,mathrsfs}

\newcommand{\n}{\nonumber}

\renewcommand{\hat}{\widehat}

\newcommand{\rbra}[1]{\!\left( #1 \right)} 
\newcommand{\cbra}[1]{\!\left\{ #1 \right\}} 
\newcommand{\sbra}[1]{\!\left[ #1 \right]} 

\newcommand{\bD}{\ensuremath{\mathbb{D}}}
\newcommand{\bE}{\ensuremath{\mathbb{E}}}

\newcommand{\bP}{\ensuremath{\mathbb{P}}}

\newcommand{\bR}{\ensuremath{\mathbb{R}}}

\newcommand{\bT}{\ensuremath{\mathbb{T}}}

\theoremstyle{plain}
\newtheorem{Thm}{Theorem}[section]

\newtheorem{Lem}[Thm]{Lemma}

\newtheorem{Cor}[Thm]{Corollary}

\theoremstyle{definition}

\newtheorem{Def}[Thm]{Definition}
\newtheorem{Rem}[Thm]{Remark}

\newcommand{\Proof}[2][Proof]{\begin{proof}[{#1}] #2 \end{proof}}

\numberwithin{equation}{section}

\makeatletter
\renewcommand\section{\@startsection {section}{1}{\z@}%
                                   {-3.5ex \@plus -1ex \@minus -.2ex}%
                                   {2.3ex \@plus.2ex}%
                                   {\normalfont\large\bf}}
\makeatother

\makeatletter
\renewcommand\subsection{\@startsection {subsection}{1}{\z@}%
                                   {-3.5ex \@plus -1ex \@minus -.2ex}%
                                   {2.3ex \@plus.2ex}%
                                   {\normalfont\normalsize\bf}}
\makeatother

\usepackage{color}
\allowdisplaybreaks[1]

\begin{document}

\begin{center}
{\Large \bf 
Generalized scale functions of standard processes with no positive jumps
}
\end{center}
\begin{center}
Kei Noba
\end{center}

\begin{abstract}
As a generalization of scale functions of spectrally negative L\'evy processes, 
we define scale functions of general standard processes with no positive jumps. 
For this purpose, we utilize excursion measures. 
Using our new scale functions we study  
Laplace transforms of hitting times, potential measures and duality.  
\end{abstract}

\section{Introduction}\label{Intro}
We first recall the basic facts of scale functions of spectrally negative L\'evy 
processes. 
Let $(X, \bP^X_x)$ with 
$X=\{X_t : t \geq 0\}$ be a spectrally negative L\'evy process 
with {$\bP^X_x(X_0=x)=1$}, 
i.e., a L\'evy process which does not have 
positive jumps and monotone paths. 
Then there correspond to it 
the Laplace exponent $\Psi$ and the $q$-scale function $W^{(q)}$ of $X$ 
for all $q\geq 0$. 
The Laplace exponent $\Psi$ is a function from $[0, \infty)$ to $\bR$ defined by 
\begin{align}
\Psi (\lambda)= \log \bE^X_0 \sbra{ e^{\lambda X_1}}, ~~~~~~~~\lambda \geq 0. 
\end{align}
The $q$-scale function $W^{(q)}$ is a function 
which is equal to $0$ on $(-\infty , 0)$, 
{is continuous on $[0, \infty)$}, and satisfies
\begin{align}
\int_0^\infty e^{-\beta x}W^{(q)}(x)dx= \frac{1}{\Psi (\beta)-q}
,~~~~~~\beta>\Phi (q), \label{102aa}
\end{align}
where $\Phi (q)= \inf \{ \lambda > 0 : \Psi (\lambda)> q  \}$ 
(see, e.g., \cite[Section 8]{Kyp} for the details).  
The scale function is useful since 
the Laplace transform of hitting times and the $q$-potential measure 
can be characterized as follows: for $b< x < a$
\begin{align}
\bE^X_x \sbra{e^{-qT^+_a} ; T^+_a < T^-_b}
&=\frac{W^{(q)}(x-b)}{W^{(q)}(a-b)}, \label{103aa}\\
\bE^X_x \sbra{\int_0^{T^+_a \land T^-_b}e^{-qt}f(X_t)dt }
&=\int_b^a f(y)\rbra{\frac{W^{(q)}(x-b)}{W^{(q)}(a-b)}W^{(q)}(a-y) - W^{(q)}(x-y)}dy,
\label{104aa}
\end{align}
where  $T^+_a$ and $T^-_b$ denote the first hitting times of 
$[a, \infty)$ and $(-\infty , b]$, respectively.
\par
Kyprianou--Loeffen(\cite{KypLoe}) has introduced the refracted L\'evy processes and 
Noba--Yano(\cite{NobYan}) generalized their results to the generalized refracted 
L\'evy processes, where the scale functions of the processes considered were defined 
and utilized. 
In particular, Noba--Yano(\cite{NobYan}) proved that the scale functions satisfy
\begin{align}
W^{(q)}(x) = \frac{1}{n^X_0 \sbra{e^{-qT^+_x}; T^+_x < \infty }},  ~~~~~q \geq 0, ~~ 
x>0, \label{105aa}
\end{align}
where $n^X_0$ is an excursion measure away from $0$ 
subject to the normalization 
\begin{align}
n^X_0\sbra{1-e^{-qT_0}}=\frac{1}{\Phi^\prime (q)},   ~~~~q > 0 .
\end{align} 
\par
In this paper, we define generalized scale functions of $\bR$-valued standard processes by 
generalizing \eqref{105aa} and {make identities which are 
analogies of \eqref{103aa} and \eqref{104aa} of generalized scale functions. 
In addition, we study a condition of generalized scale functions for two standard processes to be in duality. }
\par
The organization of this paper is as follows. 
In Section \ref{Pre}, we prepare some notation and recall preliminary facts about 
standard processes, local times and excursion measures. 
In Section \ref{scale functions}, we give the definition of the generalized scale functions and obtain identities which are 
analogies of \eqref{103aa} and \eqref{104aa}. 
In Section \ref{duality}, we characterize duality of standard processes in terms of 
scale functions. 
In Appendix \ref{Spe}, we discuss different definitions of scale functions of spectrally negative L\'evy processes.

\section{Preliminary}\label{Pre}
Let $\bD$ denote the set of functions 
$\omega: [0, \infty) \rightarrow \bR \cup \{ \partial \}$ which are c\`adl\`ag and satisfy 
\begin{align}
\omega(t)=\partial ,~~~~~~t \geq \zeta, 
\end{align}
where $\partial$ is an isolated point and $\zeta = \inf \{t>0 : \omega (t) = \partial \}$. 
Let $\cal{B}(\bD)$ denote the class of Borel sets of $\bD$ 
equipped with the Skorokhod topology. 
For $\omega \in \bD$, let
\begin{align}
&T^-_x (\omega):= \inf \cbra{ t>0  : \omega(t) \leq x}, \\
&T^+_x (\omega):= \inf \cbra{  t>0  : \omega (t ) \geq  x}, \\
&T_x (\omega):= \inf \cbra{t> 0: \omega (t) =x}. 
\end{align}
We sometimes write $T^-_x$, $T^+_x$, $T_x$ simply for 
$T^-_x(X)$, $T^+_x(X)$, $T_x(X)$, respectively, when 
we consider a process $(X , \bP^X_x)$.  
Let $\bT$ be an interval of $\bR$ and set $a_0 = \sup \bT$ and $b_0 =\inf \bT$. 
We assume that the process $(X, \bP^X_x)$ considered in this paper is a 
$\bT$-valued standard process with no positive jumps 
{with $\bP^X_x (X_0 = x)=1$}, 
satisfying the following conditions: 
\begin{itemize}
\item[(A1)]{$(x, y) \mapsto \bE^X_x\sbra{ e^{-T_y}}>0$ is a ${\cal{B}}(\bT) \times \cal{B}(\bT)$-measurable function.} 
\label{iv}
\item[(A2)]{$X$ has a reference measure $m$ on $\bT$, i.e. 
for $q \geq 0$ and $x\in \bT$, we have $R^{(q)}_X 1_{(\cdot)}(x) \ll m(\cdot)$, 
where 
\begin{align}
R_X^{(q)}f(x):= \bE^X_x \sbra{ \int_0^\infty e^{-qt} f(X_t) dt}
\end{align}
for non-negative measurable function $f$. 
Here and hereafter we use the notation $\int_b^a = \int_{(b, a]\cap\bR}$. 
In particular, $\int_{b-}^a=\int_{[b, a]\cap\bR}$. 
}
\end{itemize}
{By \cite[Theorem $18.4$]{GemHor}, 
there exist a family of processes ${\{L^{Z, x}\}}_{x\in\bT}$ with 
$L^{Z, x}={\cbra{L^{Z,  x}_t}}_{t\geq 0}$ for 
$x \in \bT$ which we call \textit{local times} 
such that the following conditions hold: 
for all $q>0$, $x \in \bT$ and non-negative measurable function $f$ 
\begin{align}
&~~~\int_0^t f(X_s)ds = \int_{\bT} f(y) L^{X, y}_{t} m (dy) ,   ~~~~\text{a.s.}  \label{201f}\\
&R_X^{(q)}f(x) 
= \int_{\bT} f(y) \bE^X_x\sbra{\int_{0}^\infty e^{-qt}dL^{X, y}_{t} } m (dy). \label{202f}
\end{align}
We have the following two cases: 
\begin{itemize}
\item{{\bf Case 1.} 
If $x\in \bT$ is regular for itself, 
this $L^{Z,x}$ is the continuous local time at $x$ (\cite[pp.216]{BluGet}). 
Note that $L^{Z,x}$ has no ambiguity of multiple constant because of 
\eqref{201f} or \eqref{202f}. }
\item{{\bf Case 2.} 
If $x\in \bT$ is irregular for itself, 
we have 
\begin{align}
L^{Z, x}_t =l^Z_x \# \{0\leq s< t : Z_s=x \},~~~\text{a.s.}
\end{align}
for some constant $l^Z_x \in (0, \infty)$. 
}
\end{itemize}}
\par
If $x$ is regular for itself, the excursion measure $n^X_x$ can be defined from 
the Poisson point process (see \cite{Ito}).
Let $\eta^{X, x}$ denote the inverse local time of $L^{X, x}${, 
i.e., $\eta^{X, x}_t = \inf \{ s> 0 :L^{X, x}_s > t \}$
}. 
Then, for all $q > 0$, we have 
\begin{align}
-\log \bE^X_0\sbra{e^{-q\eta^{X, x}(1)}}
=\delta^X_x q + n^X_x\sbra{1-e^{-qT_x} } 
\end{align}
for a non-negative constant $\delta^X_x$ called the \emph{stagnancy rate}. 
We thus have 
\begin{align}
\bE^X_x \sbra{\int_{0}^\infty e^{-qt} dL^{X, x}_t} 
= \bE^X_x\sbra{\int_0^\infty e^{-q\eta^{X, x} (s)}ds} 
=\frac{1}{\delta^X_x q +n^X_x \sbra{1 - e^{-qT_x}}}. 
\label{203z} 
\end{align}
If $x$ is irregular for itself, we define $n^X_x = \frac{1}{l^X_x}\bP^{X^x}_x$ 
where $\bP^{X^x}_x$ denotes the law of $X$ started from $x$ and stopped at $x$. 
Then 
we have 
\begin{align}
&\bE^X_x \sbra{ \int_{0-}^{\infty} e^{-qt} d L^{X, x}_t }
=l^X_x \sum_{ i = 0}^{\infty} \rbra{\bE^X_x \sbra{e^{-qT_x}}}^i
=\frac{l^X_x}{\bE^X_x \sbra{ 1-e^{-qT_x}}}= \frac{1}{n^X_x \sbra{1-e^{-qT_0} }}.  
\end{align}
\begin{Rem}
Any point $x\in\bT\backslash\{a_0\}$ cannot be a holding point. 
In fact, assume $x$ {is}. 
Then $X$ leaves $x$ by jumps (see, e.g., \cite[Theorem 1 ($vi$)]{Sal}). 
But $X$ has no positive jumps, and thus $X$ {can} not exceed $x$., 
which contradicts (A1). 
\end{Rem}
\section{Generalized scale functions}\label{scale functions}
In this section, we define generalized scale functions of 
standard processes with no positive jumps and study 
these fundamental properties. 
\begin{Def}
For $q\geq 0$ and $x, y \in \bT$, we define $q$-scale function of $X$ as 
\begin{align}
W_X^{(q)}(x, y) =
\begin{cases}
\frac{1}{n^X_y \sbra{e^{-q T^+_x }; T^+_x < \infty }},  ~~~~~~~~~~   &x\geq y,
\\
0,~~~~~~~~~~~~~~~~~~~~~~~~~~~&x < y, 
\end{cases} 
\label{301aa}
\end{align} 
where $\frac{1}{\infty}=0$.
\end{Def}
\begin{Rem}\label{Rem302aa}
All $x\in\bT\backslash \{b_0\}$ is regular for {$(x, \infty)$, 
i.e., $\bP^X_x(T^+_x=0)=1$,} thanks to the 
assumptions of no positive jumps {and of (A1)}. 
When $x$ is irregular for itself, we have $W_X^{(q)}(x, x)=l^X_x$ 
by the definition of $n^X_x$. 
\end{Rem}
\begin{Rem}
Let us characterize our scale functions of diffusion processes in terms of their 
characteristics. 
Let $m$ and $s$ be two $\bR$-valued strictly increasing continuous functions 
on the interval $[0 , \infty)$ satisfying $s(0)=0$. 
Let $X$ be a $\frac{d}{dm}\frac{d}{ds}$-diffusion process with $0$ being a 
reflecting boundary. 
Note that our $n^X_0$ coincides with the excursion measure defined in 
\cite[Definition 2.1]{Yan} up to scale transformation. 
Let $\psi^{(q)}$ denote the increasing eigenfunction 
$\frac{d}{dm}\frac{d}{ds}\psi^{(q)}=q\psi^{(q)}$ such that $\frac{d}{ds}\psi^{(q)}(0)=1$. 
In other words, 
the $\psi^{(q)}$ is the unique solution of the integral equation 
\begin{align}
\psi^{(q)}(x)=s(x)+q\int_0^x(s(x)-s(y))\psi^{(q)}(y)dm(y), ~~~~~~x\in [0,\infty) .
\end{align} 
Then, by \cite[Corollary 2.4]{Yan}, for $q>0$ and $x \in (0, \infty)$, we have 
\begin{align} 
\psi^{(q)}(x)
=\frac{1}{n^X_0\sbra{e^{-qT^+_x} ; T^+_x< \infty}}, 
\end{align}
which shows that $W_X^{(q)}(x, 0)=\psi^{(q)}(x)$. 
In particular, we have $W_X^{(0)}(x, 0)=s(x)$. 
\end{Rem}
We fix $ b, a \in \bT$ with $b<a$. 
\begin{Thm}\label{Lem102}
For $q\geq 0$ and $ x \in (b, a)$, we have 
\begin{align}
\bE^X_x \sbra{ e^{-q T^+_a}; T^+_a < T^-_b}
= \frac{W_X^{(q)}(x, b)}{W_X^{(q)}(a , b) }.\label{103}
\end{align}
\end{Thm}
\Proof{
By the strong Markov property {and since $X$ has no positive jumps}, we have
\begin{align}
n^X_b\sbra{e^{-qT^+_a} ; T^+_a < \infty}
&{=n^X_b\sbra{e^{-qT^+_x} 
\bE^X_{X_{T^+_x}} \sbra{ e^{-q T^+_a}; T^+_a < T^-_b} ; T^+_x < \infty} }
\\
&=n^X_b\sbra{e^{-qT^+_x} ; T^+_x < \infty}
\bE^X_x \sbra{ e^{-q T^+_a} ; T^+_a < T^-_b}
\end{align}
and we obtain \eqref{103}. }
	\begin{Lem}\label{Lem103}
	For $q \geq 0$ and $x\in (b, a)$, we have
		\begin{align}
		\bE^X_{x}\sbra{ e^{-qT_{a}^{+}} ;T_{a}^{+} <T_{b}^{-}}
		&=n^X_x\sbra{ e^{-qT_{a}^{+}}; T^+_a < \infty } 
\bE^X_x \sbra{\int_{0-}^{T^+_a \land T^-_b}e^{-qt}dL^{X, x}_t}. \label{211v} 
		\end{align}
	\end{Lem}
\Proof{The proof is almost the same as that of {\cite[Lemma $6.1$]{NobYan}}, 
but a slight difference lies in presence of stagnancy. \par
$i)$ We assume that $x$ is regular for itself. \par
Let $p=\{p{(t)}: t\in D(p)\}$ denote a Poisson point process 
with characteristic measure $n^X_x$. 
We extend $p(t)$ for all $t \geq 0$ by setting 
$p(t)=\partial$ for $t \notin D(p)$ 
where we abuse $\partial$ for the path taking values in $\partial$ for all time.
For $E \in {\cal{B}}(\bD)$, we write $\kappa_E = \inf \{s \geq 0 , p(s) \in E\}$. 
We let $A = \{T^+_a < \infty\} \cup \{T^-_b < \infty \}\cup\{ \zeta< \infty\}$. 
We denote by $\epsilon^\star= p(\kappa_A)$ the first excursion belonging to $A$. 
Then we have
\begin{align}
\bE^X_x \sbra{ e^{-qT^+_a}  ;T^+_a < T^-_b }
&=\bE^X_x \sbra{e^{-q \eta^{X,x} (\kappa_A -) }e^{-q T^+_a (\epsilon^\star)}
{; T^+_a(\epsilon^\ast) < T^-_b (\epsilon^\ast)}} . 
\end{align}
{Since $X$ has no positive jumps, we have $T^+_a (\epsilon^\ast)=\infty$ 
and $e^{-qT^+_a(\epsilon^\ast)}=0$ 
on $\{T^+_a (\epsilon^\ast)> T^-_b (\epsilon^\ast)\}$. 
So we have
\begin{align}
\bE^X_x \sbra{e^{-q \eta^{X,x} (\kappa_A -) }e^{-q T^+_a (\epsilon^\star)}
; T^+_a(\epsilon^\ast) < T^-_b (\epsilon^\ast)}
=\bE^X_x \sbra{e^{-q \eta^{X,x} (\kappa_A -) }e^{-q T^+_a (\epsilon^\star)}}
\end{align}
By the renewal property of the Poisson point process, we have
\begin{align}b
\bE^X_x \sbra{e^{-q \eta^{X,x} (\kappa_A -) }e^{-q T^+_a (\epsilon^\star)}}
=\bE^X_x \sbra{e^{-q \eta^{X,x} (\kappa_A -) }}
\frac{n^X_x\sbra{e^{-qT^+_a}; A}}{n^X_x \sbra{A}}.
\end{align}}
We denote ${D(p_A)}=\{s \in D(p): p(s)\in A\}$ and $p_A=p|_{D(p_A)}$ 
and extend $p_A(t)$ for all $t\geq 0$ 
{by setting $p_A(t)=\partial$ for $t \not\in D(p_A)$}. 
For $q>0$, let ${\bf{e}}_q$ be an independent random variable 
which has the exponential distribution with intensity $q$.  
We write $\eta^{X,x}_{A^c}(s)= \eta^{X,x}(s) - \sum_{u\leq s}T_x (p_A (u))$ where 
$T_x(\partial) = 0$. 
Since $\eta^{X,x} (\kappa_A -) = \eta^{X,x}_{A^c} (\kappa_A )$ where $\eta^{X,x}_{A^c}$ 
and $\kappa_A$ are independent, we have 
\begin{align}
\bE^X_x \sbra{e^{-q \eta^{X,x} (\kappa_A -) }}
&
=\bE^X_x \sbra{{\exp\rbra{-q \eta^{X,x}_{A^c} ({\bf{e}}_{n^X_x\sbra{A}}) }}}\\
&=n^X_x\sbra{A}\bE^X_x\sbra{{  \int_0^\infty 
\exp\rbra{-t n^X_x\sbra{A}}\exp\rbra{-q\eta^{X, x}_{A^c} (t)}dt }} \label{215w}\\
&=n^X_x\sbra{A}\bE^X_x\sbra{{  \int_0^{{\bf{e}}_{n^X_x\sbra{A}} } 
\exp\rbra{-q\eta^{X,x}_{A^c} (t)}dt }}\\
&=n^X_x\sbra{A}\bE^X_x\sbra{{  \int_0^{\kappa_A } 
\exp\rbra{-q\eta^{X, x} (t)}dt }} \\
&=n^X_x\sbra{A} \bE^X_x \sbra{\int_{0-}^{T^+_a \land T^-_b}e^{-qt}dL^{X, x}_t}.
\end{align}
Therefore we obtain \eqref{211v}. 
\par
$ii)$ We assume that $z$ is irregular for itself. \par 
Let $T_x^{(n)}$ denote the $n$-th hitting time to $x$ and 
let $T_x^{(0)}=0$. 
Then we have 
\begin{align}
\bE^X_x \sbra{ \int_{0-}^{T^+_a \land T^-_b} e^{-qt} dL^{X, x}_t}
&=l^X_x\sum_{i=0}^{\infty} \bE^X_x \sbra{ e^{-qT^{(i)}_x};T^{(i)}_x<T^+_a \land  T^-_b}\\
&= l^X_x \sum_{i=0}^{\infty}{\rbra{ \bE^X_x \sbra{ e^{-qT_x}; T_x<T^+_a \land  T^-_b}}}^i . 
\end{align}
On the other hand, we have 
\begin{align}
\bE^X_x \sbra{ e^{-qT^+_a};T^+_a< T^-_b }
&=\sum_{i= 0}^\infty 
{\rbra{\bE^X_x \sbra{ e^{-qT_x} ; T_x<T^+_a \land  T^-_b}}}^i
\bE^X_x \sbra{ e^{-qT^+_a} ; T^+_a < T_x \land T^-_b} . 
\end{align}
Therefore we obtain \eqref{211v}. 
}
By Theorem \ref{Lem102} and Lemma \ref{Lem103}, 
for $q \geq 0$ and $x\in (b, a)$, we obtain 
	\begin{align}
	\bE^X_x \sbra{\int_{0-}^{T^+_a \land T^-_b}e^{-qt}dL^{X, x}_t}= 
	\frac{W_X^{(q)}(x, b)W_X^{(q)}(a, x)}{W_X^{(q)}(a, b)} \label{Lem104}
	\end{align}
\par
For $q \geq 0$, $x \in (b, a)$ and non-negative measurable function $f$, we define
\begin{align}
\overline{\underline{R}}_X^{(q; b, a)}f(x)
:= \bE^X_x\sbra{ \int_0^{T^-_b \land T^+_a} e^{-qt}f(X_t)dt}.
\end{align}
Then, for $q \geq 0$, we have 
\begin{align}
\overline{\underline{R}}_X^{(q; b, a)}f(x)
&= \int_{(b, a)} f(y) \bE^X_{x} 
\sbra{\int_{0}^{T^-_b \land T^+_a} e^{-qt}dL^{X, y}_{t} } m (dy) .
\end{align}
\begin{Thm}\label{Lem105}
For $q \geq 0 $ and $x, y \in (b, a)$, we have
\begin{align}
\bE^X_{x} \sbra{\int_{0}^{T^-_b \land T^+_a} e^{-qt}dL^{X, y}_{t} }
= \frac{W_X^{(q)}(x, b)}{W_X^{(q)}(a, b)}W_X^{(q)}(a , y) - W_X^{(q)}(x, y). \label{114}
\end{align}
\end{Thm}
\Proof{
$i)$ We assume that $x=y$. 
\par
When $x$ is regular for itself, we have 
\begin{align}
\bE^X_{x} \sbra{\int_{0}^{T^-_b \land T^+_a} e^{-qt}dL^{X, y}_{t} }
=\bE^X_{x} \sbra{\int_{0-}^{T^-_b \land T^+_a} e^{-qt}dL^{X, y}_{t} }
\end{align}
and 
\begin{align}
W_X^{(q)}(x, x) = 0. 
\end{align}
By \eqref{Lem104}, we obtain \eqref{114}. 
\par
When $x$ is irregular for itself, we have 
\begin{align}
\bE^X_{x} \sbra{\int_{0}^{T^-_b \land T^+_a} e^{-qt}dL^{X, y}_{t} }
&=\bE^X_{x} \sbra{\int_{0-}^{T^-_b \land T^+_a} e^{-qt}dL^{X, y}_{t} }-l^X_x . 
\end{align}
By \eqref{Lem104} and Remark \ref{Rem302aa}, we obtain \eqref{114}. 
\par
$ii)$ We assume that $x < y$.
\par
By the strong Markov property, 
we have
\begin{align}
\bE^X_{x} \sbra{\int_{0}^{T^-_b \land T^+_a} e^{-qt}dL^{X, y}_{t} }
&=\bE^X_x\sbra{e^{-q T^+_y }; T^+_y < T^-_b}
\bE^X_{y} \sbra{\int_{0-}^{T^-_b \land T^+_a} e^{-qt}dL^{X, y}_t } . \label{213f}
\end{align}
So by Theorem \ref{Lem102} and \eqref{Lem104}, we have
\begin{align}
\eqref{213f}
=\frac{W_X^{(q)}(x ,b)}{W_X^{(q)}(y ,b)} 
\frac{W_X^{(q)}(y, b)W_X^{(q)}(a, y)}{W_X^{(q)}(a, b)}
\end{align}
and we obtain \eqref{114}. 
\par
$iii)$ We assume that $x > y$.
\par
By the strong Markov property, we have
\begin{align}
\bE^X_{x} \sbra{\int_{0}^{T^-_b \land T^+_a} e^{-qt}dL^{X, y}_{t} }
&=\bE^X_x\sbra{e^{-q T_y }; T_y < T^-_b \land T^+_a}
\bE^X_y \sbra{\int_{0-}^{T^-_b \land T^+_a} e^{-qt}dL^{X, y}_t }. \label{119}
\end{align}
Since
\begin{align}
\bE^X_x\sbra{ e^{-qT^+_a};T_y < T^+_a <T^-_b}
=\bE^X_x\sbra{e^{-q T_y }; T_y < T^-_b \land T^+_a}
\bE^X_y\sbra{ e^{-qT^+_a}; T^+_a <T^-_b} , 
\end{align} 
we have
\begin{align}
\eqref{119}
&=\frac{\bE^X_x\sbra{ e^{-qT^+_a};T_y < T^+_a <T^-_b}}
{\bE^X_y\sbra{ e^{-qT^+_a}; T^+_a <T^-_b}}
\bE^X_y \sbra{\int_{0-}^{T^-_b \land T^+_a} e^{-qt}dL^{X, y}_t }
\\
&=\frac{\bE^X_x\sbra{ e^{-qT^+_a}; T^+_a <T^-_b}
-\bE^X_x\sbra{ e^{-qT^+_a}; T^+_a <T^-_y}}
{\bE^X_y\sbra{ e^{-qT^+_a}; T^+_a <T^-_b}}
\bE^X_y \sbra{\int_{0-}^{T^-_b \land T^+_a} e^{-qt}dL^{X, y}_t }.
\label{122}
\end{align}
By Theorem \ref{Lem102} and \eqref{Lem104}, we have 
\begin{align}
\eqref{122}
=\frac{W_X^{(q)}(a, b)}{W_X^{(q)}(y, b)}
\rbra{\frac{W_X^{(q)}(x, b)}{W_X^{(q)}(a, b)}-\frac{W_X^{(q)}(x, y)}{W_X^{(q)}(a, y)}}
\frac{W_X^{(q)}(y, b)W_X^{(q)}(a, y)}{W_X^{(q)}(a, b)}
\end{align}
and we obtain \eqref{114}. 
}
\begin{Cor}
For $x, y \in (b_0, a_0)$, we define 
\begin{align}
Z_X^{(q)} (x, y)= 
\begin{cases}
1 + q\int_{(y,x)} W_X^{(q)}(x, z)m (dz),  ~~~&x > y, \\
1,~~~~~~~~~~~~~~~~~~~~~~~~~~~~~~~~~~~ &x \leq y. 
\end{cases}
\end{align}
Then we have 
\begin{align}
\bE^X_x \sbra{e^{-q T^-_b} ;T^-_b < T^+_a }
=Z_X^{(q)}(x, b) - \frac{W_X^{(q)}(x, b)}{W_X^{(q)}(a, b)}Z_X^{(q)}(a, b). \label{226i}
\end{align}
\end{Cor}
\Proof{
We have 
\begin{align}
\bE^X_x \sbra{e^{-q T^-_b} ;T^-_b < T^+_a }
=\bE^X_x \sbra{e^{-q (T^-_b \land T^+_a )} ;T^-_b \land T^+_a < \infty } 
-\bE^X_x \sbra{e^{-q T^+_a} ;T^+_a < T^-_b }. \label{225h}
\end{align}
By Theorem \ref{Lem105}, we have 
\begin{align}
&\bE^X_x \sbra{e^{-q (T^-_b \land T^+_a )} ;T^-_b \land T^+_a< \infty } 
\notag \\&
=1 - q\bE^X_x\sbra{\int_0^{T^-_b \land T^+_a}e^{-qt}dt }\\
&=1 - q\int_{(b, a)} \rbra{\frac{W_X^{(q)}(x, b)}{W_X^{(q)}(a, b)}W_X^{(q)}(a , y) - W_X^{(q)}(x, y)}
	m( dy ). \label{336ab}
\end{align}
By \eqref{225h}, \eqref{336ab} and Theorem \ref{Lem102}, we have 
\begin{align}
\eqref{225h}
=1+ q\int_{(b, a)} W_X^{(q)}(x, y) m (dy )
-\frac{W_X^{(q)}(x, b)}{W_X^{(q)}(a, b)} \rbra{1+ q\int_{(b, a)} W_X^{(q)}(a, y) m (dy ) }, 
\end{align}
and therefore we obtain \eqref{226i}. 
}

\section{Relation between duality and scale functions}\label{duality}
In this section, we characterize duality of standard processes in terms of 
our scale functions. 
\par
Let $X$ be a $\bT$-valued standard process with no positive jumps satisfying 
(A1) and (A2). 
Let $(\hat{X} , \bP^{\hat{X}}_x)$ with $\hat{X}=\cbra{\hat{X}_t: t \geq 0}$ be a 
$\bT$-valued standard process with no negative jumps 
satisfying the following conditions: 
\begin{enumerate}
\item[(B1)]{$(x, y) \mapsto \bE^{\hat{X}}_x\sbra{ e^{-T_y}}>0$ is a ${\cal{B}}(\bT) \times \cal{B}(\bT)$-measurable function.} 
\item[(B2)]{$\hat{X}$ has a reference measure $m$ on $\bT$.}
\end{enumerate}
For $q\geq 0$ and non-negative measurable function $f$, we denote 
\begin{align}
R^{(q)}_{\hat{X}} f (x)= \bE^{\hat{X}}_x \sbra{ \int_0^\infty e^{-qt}f(\hat{X}_t)dt}. 
\end{align}
We define local times ${\{L^{\hat{X}, x}\}}_{x\in\bT}$, 
excursion measures ${\{n^{\hat{X}}_x\}}_{x \in \bT}$ 
and scale functions ${\{ W_{-\hat{X}}^{(q)}\}}_{q \geq 0}$ 
of $\hat{X}$ by the same way as $X$'s in Section \ref{scale functions}.
\begin{Def}[See, e.g., \cite{Chu}]
Let $m$ be a $\sigma$-finite Radon measure on $\bT$. 
We say that $X$ and $\hat{X}$ are \emph{in duality (relative to $m$)} if 
for $q >0$, non-negative measurable functions $f$ and $g$, 
\begin{align}
\int_\bT f(x)R^{(q)}_Xg(x)m(dx)=\int_\bT R^{(q)}_{\hat{X}}f(x) g(x)m(dx). 
\end{align}
\end{Def}
\begin{Thm}[See, e.g., \cite{Chu} or \cite{Rev}]\label{Thm102a}
Suppose $X$ and $\hat{X}$ be 
in duality relative to $m$. 
Then, for each $q > 0$, there exists a function $r_X^{(q)}: \bT \times \bT \rightarrow 
[0, \infty)$ such that
\begin{enumerate}
\item{$r_X^{(q)}$ is ${\cal{B}}(\bT)\times {\cal{B}}(\bT)$-measurable.}
\item{$x \mapsto r_X^{(q)}(x, y) $ is $q$-excessive 
and finely continuous for each $y \in \bT$.}
\item{$y \mapsto r_X^{(q)}(x, y) $ is $q$-coexcessive 
and cofinely continuous for each $x \in \bT$.}
\item{For all non-negative function $f$, 
\begin{align}
R_X^{(q)}f(x)=\int_\bT f(y) r_X^{(q)}(x, y) m(dy), ~~~~~  
R_{\hat{X}}^{(q)}f(y)=\int_\bT f(x) r_X^{(q)}(x, y) m(dx). \label{402x}
\end{align}
}
\end{enumerate}
\end{Thm}
\begin{Rem}
We suppose that $X$ and $\hat{X}$ are in duality. 
Then if a function $f$ is finely continuous at $x\in \bT$, 
then $f$ is right continuous at $x$. 
If $f$ is cofinely continuous at $x\in \bT$, then $f$ is left continuous at $x$.  
\end{Rem}  
\begin{Rem}
We suppose that $X$ and $\hat{X}$ are in duality relative to $m$. 
Then, by the proof of Lemma \ref{Lem404x} and Theorem \ref{Lem105}, 
for $x, y\in (b_0 , a_0)$, the function $x\mapsto W^{(q)}(x , y)$ is finely continuous and 
$y\mapsto W^{(q)}(x , y)$ is cofinely continuous. 
\end{Rem}
By 
\cite[Proposition of Section $V$.1]{Rev}, 
if $X$ and $\hat{X}$ are in duality relative to $m$, 
there exist local times ${\{L^{X, x}\}}_{x \in \bT} $ of $X$ and 
${\{{{L}}^{\hat{X}, x}\}}_{x \in \bT} $ of $\hat{X}$ satisfying 
\begin{align}
\bE^X_x\sbra{\int_{0}^\infty e^{-qt} dL^{X, y}_{t}}=r_X^{(q)}(x , y), ~
\bE^{\hat{X}}_y\sbra{\int_{0}^\infty e^{-qt} dL^{{\hat{X}}, x}_{t}}=r_X^{(q)}(x , y). \label{403x}
\end{align}
for all $q >0$. 
When $X$ and $\hat{X}$ are in duality, we always use the normalization of the local times above. 
In other cases, we use the normalization of the local times in Section \ref{scale functions}. 
\begin{Thm}\label{Lem206b}
If $X$ and $\hat{X}$ are in duality relative to $m$, 
then we have
\begin{align}
W_X^{(q)}(x, y) = W_{-{\hat{X}}}^{(q)}(-y, -x),~~~~~x, y \in (b_0 , a_0 ) .   \label{123}
\end{align} 
If $\bT$ is open, then the converse is also true. 
\end{Thm}
\begin{Rem}
If $X$ is a spectrally negative L\'evy process, 
we have that $\hat{X}=-X$ is a dual process relative to $m$ when $m$ is the 
Lebesgue measure. 
By stationary independent increments and the definition of 
generalized scale functions, we have
\begin{align}
W_X^{(q)}(x, y)= W^{(q)}(x- y)=W_{-\hat{X}}^{(q)}(-y , -x) , 
\end{align}
where $W^{(q)}$ is the right continuous function defined by \eqref{102aa}. 
We thus obtain \eqref{123}. 
\end{Rem}
\begin{Lem}\label{Lem404x} 
Suppose $X$ and $\hat{X}$ be in duality relative to $m$. 
Then, for all $b< a \in \bT$ and $x, y\in (b , a)$, we have 
\begin{align}
\bE^X_{x} \sbra{\int_{0}^{T^-_b \land T^+_a} e^{-qt} dL^{X, y}_{t} }
&=\bE^{\hat{X}}_{y} \sbra{\int_{0}^{T^-_b \land T^+_a} e^{-qt} dL^{\hat{X}, x}_{t}}. \label{406aa}
\end{align}
\end{Lem}
\Proof{
Let $X^{(b, a)}$ and ${\hat{X}}^{(b , a)}$ denote the $X$ and $\hat{X}$ killed on exiting 
$(b , a)$, respectively. 
We denote by $R^{(q)}_{X^{(b, a)}}$ and $R^{(q)}_{{\hat{X}}^{(b, a)}}$ 
the $q$-resolvent operators of $X^{(b, a)}$ and ${\hat{X}}^{(b, a)}$, respectively. 
For each $q > 0$, there exists a function $r_{X^{(b, a)}}^{(q)}: 
(b, a) \times (b , a) \rightarrow [0, \infty)$ such that
\begin{enumerate}
\item{$r_{X^{(b, a)}}^{(q)}$ is ${\cal{B}}(b, a)\times {\cal{B}}(b ,a)$-measurable.}
\item{$x \mapsto r_{X^{(b, a)}}^{(q)}(x, y) $ is $q$-excessive 
and finely continuous 
for each $y \in (b ,a)$.}
\item{$y \mapsto r_{X^{(b, a)}}^{(q)}(x, y) $ is $q$-coexcessive 
and cofinely continuous 
for each $x \in (b , a)$.}
\item{For all non-negative function $f$, 
\begin{align}
R_{X^{(b, a)}}^{(q)}f(x)=\int_{(b , a)} f(y) r_{X^{(b, a)}}^{(q)}(x, y) m(dy), ~~~~~  
R_{{\hat{X}}^{(b, a)}}^{(q)}f(y)=\int_{(b , a)} f(x) r_{X^{(b, a)}}^{(q)}(x, y) m(dx).
\label{407x}
\end{align}
}
\end{enumerate}
By definition, we have
\begin{align}
R_{{X}^{(b, a)}}^{(q)}f(y)= \overline{\underline{R}}_{X}^{(q; b,a)}f(y)
=\int_{(b, a)} f(y) \bE^X_x \sbra{\int_{0}^{T^-_b \land T^+_a} e^{-qt}dL^{X, y}_{t}}m(dy). 
\label{408x}
\end{align}
So, for all $x\in (b, a)$, we have 
$r_{X^{(b, a)}}^{(q)}(x, \cdot) 
= \bE^X_x \sbra{\int_{0}^{T^-_b \land T^+_a} e^{-qt}dL^{X, \cdot}_{t}}$, $m$-a.e. 
Since 
\begin{align}
\bE^X_x \sbra{\int_{0}^{T^-_b \land T^+_a} e^{-qt}dL^{X, y}_{t}}
=\bE^X_x \sbra{\int_{0}^\infty e^{-qt}dL^{X, y}_{t}}
- \bE^X_x \sbra{ e^{-q(T^-_b \land T^+_a)} 
 \bE^X_{X_{T^-_b \land T^+_a}} \sbra{ \int_{0}^\infty e^{-qt}dL^{X, y}_{t}} } 
\end{align}
and the dominated convergence theorem, the function 
$y \mapsto \bE^X_x \sbra{\int_{0}^{T^-_b \land T^+_a} e^{-qt}dL^{X, y}_{t}}$ is 
a cofinely continuous function. 
By 
cofine continuity, for all $x, y \in (b , a)$, 
we have 
\begin{align}
r_{X^{(b , a)}}^{(q)} (x, y) = 
\bE^X_x \sbra{\int_{0}^{T^-_b \land T^+_a} e^{-qt}dL^{X, y}_{t}}. \label{411ab}
\end{align}
By the same way, we have 
\begin{align}
r_{X^{(b , a)}}^{(q)} (x, y) = 
\bE^{\hat{X}}_y \sbra{\int_{0}^{T^-_b \land T^+_a} e^{-qt} dL^{\hat{X}, x}_{t}}. \label{412ab}
\end{align}
By \eqref{411ab} and \eqref{412ab}, we obtain \eqref{406aa}. 
}
\Proof[Proof of Theorem \ref{Lem206b}]{
Let us assume that $X$ and $\hat{X}$ are in duality relative to $m$. 
First, we fix $b, y, a\in\bT$ with $b<y<a$. 
By Lemma \ref{Lem404x} and 
Theorem \ref{Lem105}, for all $q \geq 0$ and $x \in (b, y)$, we have 
\begin{align}
\frac{W_X^{(q)}(x, b)}{W_X^{(q)}(a, b)}W_X^{(q)}(a , y) 
=
\frac{W_{-{\hat{X}}}^{(q)}(-y, -a)}{W_{-{\hat{X}}}^{(q)}(-b, -a)}
W_{-{\hat{X}}}^{(q)}(-b , -x) . \label{125}
\end{align}
Here there exists a function $\gamma_1 :  [0, \infty) \times \bT
\rightarrow (0, \infty)$ satisfying 
\begin{align}
W_X^{(q)}(x, b)=\gamma_1 (q, b)W_{-{\hat{X}}}^{(q)}(-b, -x) ~~~~~x \in (b , a_0). 
\label{126}
\end{align}
Second, we fix $b, x, a \in \bT$ with $b<x<a$. 
For $q \geq 0$ and $y\in (x, a)$, we have \eqref{125}. 
Thus there exists a function $\gamma_2 : [0, \infty) \times \bT
\rightarrow (0, \infty)$
\begin{align}
W_X^{(q)}(a , y)=\gamma_2 (q, a)W_{-{\hat{X}}}^{(q)}(-y, -a) ~~~~~y \in (b_0 , a). 
  \label{127}
\end{align}
By \eqref{126} and \eqref{127}, for $q\geq0$ and $a, b \in (b_0, a_0)$, 
we have $\gamma_1 (q, b) = \gamma_2 (q, a) $, so $\gamma_1$ and $\gamma_2$ depend on only $q \geq 0$. 
By \eqref{125}, $\gamma_1 = \gamma_2 \equiv 1$. 
Thus, for $y, x\in (b_0 , a_0)$ with $y<x$, 
we have $W_X^{(q)}(x, y) = W_{-\hat{X}}^{(q)}(-y, -x)$. 
By the fine continuity of $W_X^{(q)}$ and 
the cofine continuity $W_{-\hat{X}}^{(q)}$, for $x\in (b_0 , a_0)$, 
we have $W_X^{(q)}(x, x) = W_{-\hat{X}}^{(q)}(-x, -x)$. 
\par
Let us assume that $\bT$ is open and that \eqref{123} is satisfied. 
Then, for $b<a\in\bT$ and $x, y\in(b, a)$, we have
\begin{align}
&\frac{W_X^{(q)}(x, b)}{W_X^{(q)}(a, b)}W_X^{(q)}(a , y) - W_X^{(q)}(x, y) \n \\
&=\frac{W_{-{\hat{X}}}^{(q)}(-y, -a)}{W_{-{\hat{X}}}^{(q)}(-b, -a)}
W_{-{\hat{X}}}^{(q)}(-b , -x) - W_{-{\hat{X}}}^{(q)}(-y, -x).\label{417ab}
\end{align}
By Theorem \ref{Lem105}, the first term and the second term of \eqref{417ab} are 
potential densities of $X$ and $\hat{X}$ killed on exiting $(b, a)$, respectively. 
We therefore conclude the duality of killed processes, 
which yields that of the original processes. 
}

\appendix
\section{The case of spectrally negative L\'evy processes}\label{Spe}
When $X$ is a spectrally negative L\'evy process, 
the usual definition of the scale function $W^{(q)}(x)$ is based on the 
Laplace transform 
\begin{align}
\int_0^\infty e^{-qt}W^{(q)}(x) dx =  \frac{1}{\Psi (\beta ) - q},~~~~~~
\beta> \Phi (q).    \tag{\ref{102aa}}
\end{align}
It satisfies 
\begin{align}
W^{(q)}(x)=e^{\Phi (q)x}r^{(q)}(0+)-r^{(q)}(-x) \label{A01}
\end{align}
{where $r^{(q)}$ is the right-continuous potential density of $X$ 
with respect to the Lebesgue measure (see \cite{Pis})} 
and, as we have mentioned it in Section \ref{Intro}, it also satisfies 
\begin{align}
W^{(q)}(x) = \frac{1}{n^X_0 \sbra{e^{-qT^+_x}; T^+_x < \infty }}. \tag{\ref{105aa}}
\end{align}
Pistorius(\cite{Pis}) provided a potential theoretic viewpoint for the scale functions 
in the sense that he started from \eqref{A01} and proved \eqref{102aa}. 
We now provide another viewpoint. 
 \par
Let $X$ be a spectrally negative L\'evy process and $m$ are 
the Lebesgue measure. 
By Section \ref{scale functions} and \ref{duality}, we can define 
local times, excursion measures and scale functions. 
Since $X$ has stationary independent increment property, for $q\geq 0$ and 
$x, y \in \bR$, we have 
\begin{align}
\bE^X_0\sbra{\int_{0}^\infty e^{-qt}dL^{X,y-x}_{t}}
=\bE^X_x\sbra{\int_{0}^\infty e^{-qt}dL^{X,y}_{t}}.
\end{align}
So there exists a 
left continuous function $r^{(q)}:\bR \rightarrow [0 , \infty)$ and 
a right continuous function 
$W^{(q)} : \bR \rightarrow [0 , \infty)$ satisfying 
\begin{align}
r^{(q)} ( y -x)&=\bE^X_x\sbra{\int_{0}^\infty e^{-qt}dL^{X, y}_{t}}, ~~~~~x, y \in \bR, \\
W^{(q)} ( x - y)&= W_X^{(q)} (x, y), ~~~~~x , y \in \bR , 
\end{align}
for all $q \geq 0$. 
Then $r^{(q)}$ is a c\`agl\`ad function. In fact, we have 
\begin{align}
\lim_{h\downarrow 0}r^{(q)}(x+h)
&=\lim_{h\downarrow 0}\bE^X_{-h}\sbra{\int_0^\infty e^{-qt}dL^{X, x}_t}\\
&=\lim_{h\downarrow 0}\bE^X_{-h}\sbra{e^{-qT^+_0}:T^+_0\leq T_x}\bE^X_{0}\sbra{\int_{0-}^\infty e^{-qt}dL^{X, x}_t} \n\\
&~~~~~~~~+\lim_{h\downarrow 0}\bE^X_{-h}\sbra{e^{-qT_x};T_x<T^+_0}\bE^X_{x}\sbra{\int_{0-}^\infty e^{-qt}dL^{X, x}_t}\\
&=\bE^X_{0}\sbra{\int_{0-}^\infty e^{-qt}dL^{X, x}_t}\label{A07ab}
\end{align}
where \eqref{A07ab} uses 
\begin{align}
\bE^X_0\sbra{e^{-qT^+_x} ; T^+_x < \infty}=e^{-\Phi (q)x}, ~~~~~~x> 0, \label{A08}
\end{align}
and $\lim_{x\downarrow 0}\bE^X_0\sbra{e^{-qT^+_x} ; T^+_x < \infty}=\lim_{x\downarrow 0}e^{-\Phi (q)x}=1$ 
(see, e.g., \cite[Theorem 3.2]{Kyp}). 
\begin{Thm}\label{ThmA01}
For all $q \geq 0$, we have \eqref{102aa}. 
\end{Thm}
\Proof{
By the equation obtained from \eqref{211v} when $b$ limits to infinity, 
for $x>0$, 
we have 
\begin{align}
W^{(q)}(x) &=\frac{1}{n^X_0\sbra{e^{-qT^+_x} ; T^+_x < \infty} } \\
&=\frac{1}{\bE^X_0 \sbra{e^{-qT^+_x}; T^+_x <\infty}}
\bE^X_0 \sbra{\int_{0-}^{T^+_x}e^{-qt}dL^0_t}\\
&=\frac{1}{\bE^X_0 \sbra{e^{-qT^+_x}; T^+_x <\infty}}
\rbra{ \bE^X_0 \sbra{\int_{0-}^{\infty}e^{-qt}dL^0_t} -
\bE^X_0\sbra{ e^{-qT^+_x} ; T^+_x < \infty} \bE^X_x \sbra{\int_{0}^{\infty} e^{-qt}dL^0_t} } \\ 
&=\frac{1}{\bE^X_0 \sbra{e^{-qT^+_x}; T^+_x <\infty}}r^{(q)}(0+) -r^{(q)}(-x).
\end{align}
By \eqref{A08}, for $\beta > \Phi (q)$, we have 
\begin{align}
\int_0^\infty e^{-\beta x}W^{(q)}(x)dx &= 
\int_0^\infty \rbra{e^{-(\beta-\Phi (q))x} r^{(q)}(0+) -e^{-\beta x} r^{(q)}(-x)}dx\\
&=\frac{r^{(q)}( 0+)}{\beta - \Phi (q)}-\int_0^\infty e^{-\beta x}r^{(q)}(-x)dx . \label{A10}
\end{align} 
On the other hand, for all $\beta < \Phi (q)$, we have 
\begin{align}
\frac{1}{q-\Psi (\beta)}
&=\int_0^\infty e^{-(q - \Psi (\beta))t}dt \\
&=\int_0^\infty e^{-qt}\bE^X_0 \sbra{e^{\beta X_t}}dt \\
&=\int_{-\infty}^\infty e^{\beta x} r^{(q)}(x) dx \\
&=\int_0^\infty e^{\beta x}\bE^X_0 \sbra{e^{-qT^+_x};T^+_x< \infty}r^{(q)}(0+) dx
+\int_0^\infty e^{-\beta x}r^{(q)}(-x) dx \\
&=\int_0^\infty  e^{-(\Phi (q)-\beta) x}r^{(q)}(0+) dx
+\int_0^\infty e^{-\beta x}r^{(q)}(-x) dx  \\
&=\frac{r^{(q)}(0+)}{\Phi (q)-\beta} 
+\int_0^\infty e^{-\beta x}r^{(q)}(-x) dx . \label{A16}
\end{align}
By analytic extension, we have \eqref{A16} for $\beta>\Phi (q)$. 
By \eqref{A10} and \eqref{A16}, we obtain \eqref{102aa}. 
}
\begin{Rem}
The latter part of the proof of Theorem \ref{ThmA01} 
is almost the same as a part of the proof of \cite[Theorem 1($i$)--($iii$)]{Pis}.
\end{Rem}

\section*{Acknowledgments}
I would like to express my deepest gratitude to 
Professor V\'ictor Rivero, Professor Kazutoshi Yamazaki 
and my supervisor Professor Kouji Yano 
for comments and improvements. 
Especially, Professor Kouji Yano gave me a lot of advice. 
The author was supported by JSPS-MAEDI Sakura program. 


\end{document}